\magnification 1200
\def\R{{\rm I\kern-0.2em R\kern0.2em \kern-0.2em}}
\def\N{{\rm I\kern-0.2em N\kern0.2em \kern-0.2em}}
\def\P{{\rm I\kern-0.2em P\kern0.2em \kern-0.2em}}
\def\B{{\rm I\kern-0.2em B\kern0.2em \kern-0.2em}}
\def\C{{\bf \rm C}\kern-.4em {\vrule height1.4ex width.08em depth-.04ex}\;}

\def\D{{\Delta}}

\def\z{{\zeta}}
\def\cC{{\cal C}}
\def\cP{{\cal P}}

\font\ninerm=cmr8
\centerline {\bf A DECOMPOSITION OF FUNCTIONS WITH}

\centerline {\bf ZERO MEANS ON CIRCLES}
\vskip 4mm
\centerline {Josip Globevnik}
\vskip 4mm
{\noindent \ninerm ABSTRACT\  It is well known that every H\" older continuous function 
on the unit circle is the sum of two functions such that one of these functions
extends holomorphically into 
the unit disc and the other extends holomorphically into the complement of the unit disc. 
We prove that an analogue of this holds for H\" older continuous functions on an annulus  
A  which have zero averages on all circles contained in  A  which surround the hole.} 
\vskip 4mm
\bf 1.\ Introduction and the main results \rm 
\vskip 2mm
Given $a\in\C,\ \rho >0$  write $\D (a,\rho )= \{\z\in\C\colon\ |\z -a|<\rho \}$ 
and $\D = \D (0,1)$. Denote by $\B$ the open unit ball in $\C^2$. A function $f$ on 
a set $K\subset \C$ is called H\" older continuous on $K$ (with exponent $\alpha$) 
if there are
constants $M<\infty$ and $\alpha,\ 0<\alpha <1$, such that $|f(z)-f(w)|\leq M |z-w|^\alpha\ \
(z,w\in K)$. 
\vskip 1mm
Let  $a\in\C,\ \rho >0$ and let $f$ be a H\" older 
continuous function $f$ on $b\D (a,\rho )$. 
 It is well known that 
$$
f=f^+ + f^-
\eqno (1.1)
$$
where $f^+$ and $f^-$ are H\" older continuous functions on $b\D (a,\rho )$ such that 
$f^+$ 
$$
\hbox{has a continuous extension to\ }\overline\D (a,\rho)\hbox{
\ which is holomorphic on\ }\D (a,\rho )
\eqno (1.2)
$$
and $f^-$ 
$$
\eqalign{&\hbox{has a continuous extension 
to\ }[\C\cup\{\infty\}]\setminus \D (a,\rho)\hbox{
\ which is holomorphic\ }\cr
&\hbox{on\ }
[\C\cup\{\infty\}]\setminus \overline\D (a,\rho)\hbox{\ and 
vanishes at\ }\infty \cr}
\eqno (1.3)
$$
and this decomposition is unique. In fact,
$$
{1\over {2\pi i}} \int_{b\D(a,\rho )} {{f(\z )d\z }\over{\z-z}}=
\left\{\eqalign{ 
&f^+(z)\ \ \ \ \ \ \ \ (z\in\D (a,\rho)) \cr
&-f^-(z)\ \ \ \ (z\in
[\C\cup\{\infty\}]\setminus \overline\D (a,\rho).\cr}\right\}
\eqno (1.4)
$$
\vskip 2mm
In the present paper we consider functions on the annulus 
$$
A=\{\z\in\C\colon \ r_1\leq 
|\z |\leq r_2\} 
$$
where $0<r_1<r_2<\infty $ and ask for a 
decomposition similar to (1.1). Suppose that $f$ is a H\" older continuous 
function on the annulus $A$. 
For every circle $b\D (a,\rho)\subset A$ surrounding the
origin we have 
$f|b\D (a,\rho)= f^+_{a,\rho} + f^-_{a,\rho}$ where 
$f^+_{a,\rho}$ satisfies (1.2) and  $f^-_{a,\rho}$ 
satisfies (1.3).  In general 
there are no functions $f^+$ and $f^-$ on 
$A$ such that $f^+|b\D (a,\rho) = f^+_{a,\rho}$ 
and $f^-|b\D (a,\rho) = f^-_{a,\rho}$ 
whenever $b\D (a,\rho )\subset A$ surrounds the origin.
In the present paper we 
prove that there are such functions $f^+$ and $f^-$ on 
$A$ whenever $f$ satisfies 
$$
{1\over {2\pi}}\int_0^{2\pi} f(a+\rho e^{i\theta })d\theta = 0
\eqno (1.5)
$$
for every $b\D (a,\rho )\subset A$ which surrounds the origin:    
\vskip 2mm
\noindent \bf Theorem 1.1\ \ \it Let $f$ be a H\" older continuous 
function on $A$ which satisfies (1.5) whenever $b\D (a,\rho )\subset A$ 
surrounds the origin. Then $f=f^++f^-$ where $f^+$ and $f^-$ are H\" older 
continuous functions on $A$ such that for each $b\D (a,\rho )\subset A$ 
which surrounds the origin, $f^+|b\D (a,\rho)$ satisfies (1.2) 
and $f^-|b\D (a,\rho)$ 
satisfies (1.3). \rm
\vskip 1mm
We will also show that, as in the case of the circle, 
we can view $f^+$ and $f^-$ as the 
boundary values of functions, holomorphic on appropriate 
domains. To describe this,
we first rewrite the circle case in a form suitable for 
generalization. 

Let $f$ be a continuous function on $b\D (a,\rho)$. The
idea is to define a new function $F$ 
on $\{(\z,\overline\z)\colon\ \z\in b\D(a,\rho)\}$, that
is, on $b\D (a,\rho)$  "lifted" 
to 
$$
\Sigma =\{ (\z,\overline\z)\colon\ \z\in \C\},
$$ by 
$$
F(\z ,\overline\z)= f(\z )\ \ \ \ (\z\in b\D (a,\rho ))
$$
[G2] and then to write $F$ as the sum of boundary values of
holomorphic functions.  

Let
$$
\Lambda _{a,\rho} = \{ (z,w)\in\C ^2\colon\ (z-a)(w-\overline a)=\rho ^2\}.
$$
The intersection $\Lambda _{a,\rho}\cap\Sigma$ is the circle
$\{ (\z,\overline\z)\colon\ \z\in b\D (a,\rho )\}$  
whose complement in $\Lambda _{a,\rho}$ has
two components, $\Lambda _{a,\rho}^+$ and 
$\Lambda _{a,\rho}^-$ where 
$$
\eqalign{
\Lambda _{a,\rho}^+ &= 
\{(z,w)\colon\ (z-a)(w-\overline a)=\rho ^2,\ 0<|z-a|<\rho)\} \cr
\Lambda _{a,\rho}^- &= 
\{(z,w)\colon\ (z-a)(w-\overline a)=\rho ^2,\ \rho<|z-a|)\} \cr
                    &= 
					\{(z,w)\colon\ (z-a)(w-\overline a)=
					\rho ^2,\  0<|w-a|<\rho)\} \cr
					&=
					\{ (z,w)\colon\ (\overline w,\overline z)
					\in\Lambda _{a,\rho}^+\}.\cr}
$$
The sets $\Lambda _{a,\rho }^+$ and $\Lambda _{a,\rho }^-$ are 
closed one
dimensional complex submanifolds of $\C^2\setminus\Sigma$ 
attached to $\Sigma $ along
$$
b\Lambda _{a,\rho}^+ =b\Lambda _{a,\rho}^- =
\{(\z,\overline\z)\colon\
\z\in b\D (a,\rho)\}.
$$
It is easy to see that a continuous function $h$ 
on $b\D (a,\rho )$ satisfies (1.2) if 
and only 
if the function $H$ defined on $b\Lambda _{a,\rho}^+ = b\Lambda _{a,\rho}^-$ 
by $H(\z,\overline\z)=
h(\z ) \ (\z\in b\D (a,\rho))$\ has a 
bounded continuous extension from $b\Lambda _{a,\rho}^+$ 
to $\Lambda _{a,\rho}^+\cup b\Lambda
_{a,\rho}^+$ which is holomorphic on $\Lambda _{a,\rho}^+$ [G2]. 
Similarly, $h$ satisfies (1.3) 
if and only if $H$ has  has a 
bounded continuous extension from $b\Lambda _{a,\rho}^-$ to 
$\Lambda _{a,\rho}^-\cup b\Lambda
_{a,\rho}^-$ which is holomorphic on $\Lambda _{a,\rho}^-$ 
and vanishes at $\infty$.

We now pass to functions on $A$ which we will view 
as functions on 
$$
\tilde A = \{ (\z,\overline\z )\colon\ \z\in A\}.
$$
The set $\tilde A\subset\Sigma $ will be a common part of 
the boundaries of two domains 
$\Omega ^+ (A)$ and $\Omega ^- (A)$ which we now describe. 
Let $\Omega ^+(A)$ be the union of all 
$\Lambda _{a,\rho}^+$ such that $b\D (a,\rho )\subset\hbox{Int}A$ 
surrounds the origin. 
Similarly, let $\Omega ^-(A)$ be the union of all 
$\Lambda _{a,\rho}^-$ such that $b\D (a,\rho )\subset\hbox{Int}A$ 
surrounds the 
origin. Clearly, $\Omega ^- (A) $ is the image of $\Omega ^+ (A) $
under the 
reflection $(z,w)\mapsto (\overline w, \overline z)$. It turns out 
that $\Omega^+ (A)$ and
$\Omega ^- (A)$ are disjoint domains in $\C^2\setminus\Sigma$ 
attached to $\Sigma $ along $\tilde A$. For each $\z\in\hbox{Int}A$ 
there are a neighbourhood $U\subset\Sigma$ 
of $(\z,\tilde \z)$ and 
a wedge with the edge $U$ which is contained in $\Omega ^+ (A) $. 
An analogous
statement holds for $\Omega^- (A)$. 
\vskip 2mm
\noindent \bf Theorem 1.2\ \it Let $f$ be a H\" older continuous 
function on $A$ which satisfies (1.5) whenever $b\D (a,\rho )\subset A$ 
surrounds the origin.  
There are a bounded continuous function $G^+$ on 
$\Omega ^+(A)\cup b\Omega ^+ (A)$ which is holomorphic 
on $\Omega ^+(A)$ and 
a bounded continuous function $G^-$ on 
$\Omega ^-(A)\cup b\Omega ^- (A)$ which is holomorphic 
on $\Omega ^-(A)$ such that 
$$
f(z) = (1/\overline z)G^+(z,\overline z) + 
(1/z)G^- (z,\overline  z)\ \ (z\in A).
$$
\vskip 2mm \rm 
\noindent Thus, on $\tilde A$ the function $F(z,\overline z)=
f(z)$ is the sum of the boundary values of 
holomorphic functions $(1/w)G^+(z,w)$ and $(1/z)G^-(z,w)$.
\vskip 4mm
\bf 2.\ Fourier coefficients of functions with zero means \rm 
\vskip 2mm
Suppose that $f$ is a continuous function on $A$. For each
$r,\ r_1\leq r\leq r_2$, let 
$$
c_k(r)= {1\over {2\pi}}\int_{-\pi}^{\pi} 
e^{-ik\theta}f(re^{i\theta})d\theta \ \ (k\in Z)
$$
so that $\sum_{k=-\infty}^\infty c_k(r)e^{ik\theta}$  
is the Fourier series of the function $
e^{i\theta }\mapsto f(re^{i\theta })$.

We shall need the following description of the Fourier 
coefficients of functions with zero means on circles.  
\vskip 2mm
\noindent\bf Theorem 2.1\ \rm [G1, EK, V]\ \ \it A continuous 
function $f$ on $A$ satisfies (1.5) for each $b\D (a,\rho )$ $
\subset A$ surrounding the origin if and only if 

(a)\ $c_0(r)=0 \ \ \ (r_1\leq r\leq r_2)$

(b)\ for each $n\in Z,\ n\not= 0$, there are numbers 
$a_{n,0},\ a_{n,1},\ \cdots , a_{n, |n|-1}$ such that $c_n(r)=r^{-|n|}\bigl(a_{n,0}+
a_{n,1}r^2+\cdots\  a_{n, |n|-1}r^{2(|n|-1)}\bigr)\ \ (r_1\leq r\leq r_2)$.
\vskip 2mm \rm 
\noindent In the rest of this section we assume that 
$f$ is a continuous function on $A$ which satisfies 
(1.5) for each $b\D (a,\rho )\subset A$ surrounding the origin. 

If $n\geq 1$ then writing $z=re^{i\theta }$ we get 
$$\eqalign{
c_n(r)e^{i n\theta } &= r^{-n}c_n(r)z^n \cr
 &= (a_{n,0}r^{-2n}+a_{n,1}r^{-2(n-1)}+\cdots + 
 a_{n,n-1}r^{-2})z^n \cr
 &= \biggl( {{a_{n,0}}\over {z^n\overline z^n}}+
 {{a_{n,1}}\over {z^{n-1}\overline z^{n-1}}}
+\cdots + {{a_{n,n-1}}\over{z\overline z}}\biggr)z^n\cr 
&= {1\over{\overline z}}\biggl(a_{n,0}
{1\over{\overline z^{n-1}}}+ a_{n,1}{1\over{\overline z{
^{n-2}}}}z + \cdots + a_{n,n-1}z^{n-1}\biggr)\cr
&= (1/\overline z) P_{n-1}(z,1/\overline z)\cr}
$$
where $P_{n-1}$ is a homogeneous polynomial of degree $n-1$. 

If $n\leq -1$\  then we get 
$$\eqalign{
c_n(r)e^{i n\theta } &= (a_{n,0}+a_{n,1}r^2+\cdots + 
a_{n,n-1}r^{2(-n-1)})z^{-|n|} \cr 
&= {1\over{z}}\biggl(a_{n,0}{1\over{z^{|n|-1}}}+ a_{n,1}{1\over{z
^{|n|-2}}}\overline z + \cdots + a_{n,n-1}\overline z^{n-1}\biggr)\cr
&= (1/z) Q_{|n|-1}(\overline z,1/z)\cr}
$$
where $Q_{|n|-1}$ is a homogeneous polynomial of degree $|n|-1$. 
Thus, putting $z=re^{i\theta }$ into the series 
$$
(1/\overline z) \sum_{n=0}^\infty P_n(z,1/\overline z)+ (1/z)
\sum _{n=0}^\infty Q_n(\overline z,1/z)
\eqno (2.1)
$$
we get the Fourier series of the function $e^{i\theta}\mapsto f(re^{i\theta}), \ 
r_1\leq r\leq r_2$. 

For each $t,\ 0<t<1$, define the functions $f_t^+$ and $f_t^-$ on $A$ as follows
$$
f_t^+(re^{i\theta })= \sum _{k=1}^\infty t^kc_k(r)e^{ik\theta} = (1/\overline z) \sum
_{j=0}^\infty t^{j+1}P_j(z,1/\overline z)
\eqno (2.2)
$$
$$
f_t^-(re^{i\theta })= \sum _{k=-\infty}^{-1} t^{|k|}c_k(r)e^{ik\theta} = (1/z) \sum
_{j=0}^\infty t^{j+1}Q_j(\overline z,1/z)
\eqno (2.3)
$$
where $z=re^{i\theta}\in A$ \ and let
$$
f_t(z) = f_t^+ (z) + f_t^- (z)\ \ (z\in A).
$$
Note that for each $t,\ 0<t<1$, both series above converge uniformly on $A$. Write 
$$
\Phi _r(z) = {1\over{2\pi i}}\int_{b\D (r,0)}{{f(\z )d\z}\over
{\z - z}}\ \ \ (r_1\leq r\leq r_2) 
$$
and observe that 
$$
\left.\eqalign{\Phi _r (tre^{i\theta }) = f_t^+(re^{i\theta})\ \ \ 
(r_1\leq r\leq r_2,\ 0<t<1,\ \theta \in \R )\cr
\Phi _r ((1/t)re^{i\theta }) = -f_t^-(re^{i\theta})\ \ \ 
(r_1\leq r\leq r_2,\ 0<t<1,\ \theta \in \R )\cr}\ \ \ \ \ \ \right\}
\eqno(2.4)
$$
\vskip 4mm
\bf 3.\ Proof of Theorem 1.1 \rm
\vskip 2mm
Let $f$ be a H\" older continuous function on $A$ which satisfies (1.5) for each 
$b\D (a,\rho )$ which surrounds the origin. Then there are $P_n $ and $Q_n$ as above such that 
(2.1) 
with $z=re^{i\theta}$ is the Fourier series of $e^{i\theta }\mapsto f(re^{i\theta })$,
\ $r_1\leq r\leq r_2$. Using 
the decomposition (1.1) on each circle 
$b\D (0,r),\ r_1\leq r\leq r_2$, we can write $f= f^++f^-$ 
where for each $r,\ r_1\leq r\leq r_2$,
$f^+|b\D (0,r)$ satisfies (1.2) and $f^-|b\D (0,r)$ satisfies (1.3). In fact,
the functions $f^+|b\D (0,r)$ and $-f^-|b\D (0,r)$ are the limiting values of 
$\Phi _r(z)$
as $|z|\nearrow r$ or $|z|\searrow r$, respectively. Since 
$f$ is H\" older continuous on $A$
it follows that $f^+$ and $f^-$ are H\" older continuous on $A$\ [M,\ Sections 19, 20]. 

For each $r,\ r_1\leq r\leq r_2$, the function $\Phi _r$ is the Cauchy integral 
and hence 
the Poisson integral of $f^+|b\D (0,r)$ so (2.4) implies that for each 
$r,\ r_1\leq r\leq r_2$,\ \ $f_t^+(re^{i\theta})$ converges uniformly in 
$\theta$ to $f^+(re^{i\theta})$ as $t\nearrow 1$, and since $f$ is uniformly continuous on 
$A$, the standard proof of the boundary continuity of the Poisson integral shows that 
the convergence is uniform also in $r,\ r_1\leq r\leq r_2$. So $f_t$ 
converges to $f^+$ uniformly on $A$ as $t\nearrow 1$. 
Similarly we show that $f_t^-$ 
converges to $f^-$ uniformly on $A$ as $t\nearrow 1$. 

Observe that for each $b\D (a,\rho )$ that surrounds the origin,
the restriction of $1/\overline z$ to $b\D (a,\rho)$ satisfies (1.2) and the restriction of 
$\overline z$ to $b\D (a,\rho )$ has a continuous extension to $[\C\cup\{\infty\}]\setminus
\D$ which is holomorphic on $[\C\cup\{\infty\}]\setminus\overline\D$. The uniform convergence of the series
(2.2) implies that for each $t,\ 0<t<1$, $f_t^+|b\D (a,\rho)$ satisfies (1.2) whenever 
$b\D (a,\rho )\subset A$ surrounds the origin. Similarly, the uniform convergence of
the series
(2.3) and the multiplication with $1/z$ imply
that for each $t,\ 0<t<1$, $f_t^-|b\D (a,\rho)$ satisfies (1.3) whenever 
$b\D (a,\rho )\subset A$ surrounds the origin. The uniform convergence
of $f_t^+$ to $f^+$ and $f_t^-$ to $f^-$ as $t\nearrow 1$ imply that the 
analogous statements hold for $f^+$ and $f^-$. This completes the proof. 
\vskip 2mm
\noindent \bf Remark 1.\ \rm Note that the proof of
Theorem 1.1. becomes  simpler in the special case when $f$ 
is smooth, say of class $\cC ^2$ on $A$. Recall that
putting $z=re^{i\theta }$ into the series (2.1) 
we get the Fourier series of the function 
$e^{i\theta }\mapsto f(re^{i\theta })$.  Integrating by 
parts we see  
that there is a constant $M<\infty$ such that
each of the series in (2.1) is dominated on 
$A$ by the series $\sum_{n=1}^\infty Mn^{-2}$ 
which implies that one can define 
$$
f^+ (z)= (1/\overline z) \sum _{n=0}^\infty P_n(z,1/\overline z),\ \ 
f^- (z)= (1/z) \sum _{n=0}^\infty Q_n(\overline z,1/z)\ \ (z\in A)
$$
where each of the series converges uniformly on $A$. 
\vskip 2mm
\noindent \bf Remark 2.\ \rm If $f$ in Theorem 1.1 is H\" older continuous on 
$A$ with exponent $\alpha ,\ 0<\alpha<1$, then for any $\beta ,\ 0<\beta<\alpha $, 
the functions $f^+$ and $f^-$ are H\" older continuous on $A$ with exponent $\beta$. 
This follows from [M, Sections 19, 20]. 
\vskip 4mm
\bf 4.\ Domains $\Omega ^+(A)$ and $\Omega ^-(A)$\rm
\vskip 2mm
We list some simple facts about the domain $\Omega ^+(A)$. Analogous statements hold for 
$\Omega ^-(A)$, the image of $\Omega ^+(A)$ under the reflection  $(z,w)\mapsto 
(\overline w, \overline z)$. The proofs are elementary, they can be found in [G2]. 

Recall that $\Omega ^+(A)$ is defined as the union of all $\Lambda _{a,\rho}^+$ such that 
$b\D (a,\rho)\subset\hbox{Int}A $ surrounds the origin. 
\vskip 2mm
\noindent \bf Proposition 4.1\ \it Let $\gamma = (r_1+r_2)/2$. The set $\Omega ^+(A)$ 
is a disjoint union of all $\Lambda _{a,\gamma }^+$ such that
$b\D (a,\gamma)\subset\hbox{Int}A$; 
it is an unbounded open connected set whose boundary consists of $\tilde A$ together 
with the union of all those $\Lambda _{a,\gamma }^+$ for which $b\D (a,\gamma )\subset A$ is 
tangent to both $b\D (0,r_1)$ and $b\D (0, r_2)$. For each $\z\in\hbox{Int}A$ there are a 
neighbourhood $U\subset\Sigma$ of $(\z,\overline \z)$ an open cone $V$ in $i\Sigma$, a real
two-plane perpendicular to $\Sigma$, and a $\delta >0$ such that
$U+(V\cap\delta \B) \subset \Omega ^+(A)$. 
\vskip 2mm
\noindent\bf  Proposition 4.2\ \it If  $\Lambda _{a,\rho}^+\not=
\Lambda _{b,\delta}^+$ then the sets $\Lambda _{a,\rho}^+$ 
and $\Lambda _{b,\delta}^+$ intersect if and only if 
$a\not= b$  and one of the circles $b\D (a,\rho ),\ b\D (b,\delta )$ surrounds the other. 
The sets $\Lambda _{a,\rho}^+$ and $\Lambda _{b,\delta}^-$ intersect if and only if 
$\overline\D (a,\rho)\cap\overline \D (b,\delta )=\emptyset $. \rm
\vskip 2mm
\noindent Proposition 4.2 implies that $\Omega ^+(A)\cap \Omega ^-(A)=\emptyset$. 
\vskip 4mm
\bf 5.\ Proof of Theorem 1.2 \rm
\vskip 2mm
Recall that $f^+$ is the uniform limit of $f_t^+$ as $t\nearrow 1$ where for
each $t,\ 0<t<1$, 
$$
f_t^+(z)=(1/\overline z)\sum _{j=0}^\infty t^{j+1}P_j(z,1/\overline z)\ \ (z\in A)
$$
with the series converging uniformly on $A$. It follows that on $A,\ \overline zf^+(z)$ 
is a uniform limit of a sequence of polynomials in $z$ and $1/\overline z$, 
$$
\overline zf^+(z) = \lim_{m\rightarrow\infty}S_m(z,1/\overline z)\ \ (z\in A).
\eqno (5.1)
$$
Now we reason as in [G2]:\ The functions $S_m(z,1/w)$ are bounded and continuous on
$\Lambda _{a,\rho}^+\cup b\Lambda _{a,\rho}^+$ and holomorphic on $\Lambda _{a,\rho}^+$ 
whenever $b\D (a,\rho )\subset A$ surrounds the origin. Since $\Lambda _{a,\rho}^+$ 
is biholomorphically equivalent to the punctured disc 
the maximum principle implies that for each $(z,w)\in\Lambda _{a,\rho}^+$ we have 
$$
\eqalign{
|S_m(z,1/w)-S_j(z,1/w)|&\leq \max \{|S_m(\z,1/\eta)-S_j(\z, 1/\eta )|\colon\ 
(\z,\eta)\in b\Lambda _{a,\rho}^+\} \cr
&\leq \max \{\{|S_m(\z,1/\overline\z)-S_j(\z, 1/\overline\z )|\colon\ \z\in A\}.\cr}
$$
It follows that the sequence $S_m(z,1/w)$ converges 
uniformly on $\Omega ^+(A)\cup b\Omega ^+(A)$ 
to a function $G^+(z,w)$. Since each $S_m(z,1/w)$ is bounded and continuous on
$\Omega ^+ (A)\cup b\Omega ^+(A)$ and holomorphic on $\Omega ^+ (A)$ the same is true for 
$G^+$. Obviously, $f^+(z)=(1/\overline z)G^+(z,\overline z)\  (z\in A)$.
In the same way we prove that $f^-(z) =
(1/z)G^- (z,\overline z)\ (z\in A)$ where $G^-$ is bounded and continuous 
on $\Omega ^-(A)\cup b\Omega^-(A)$ and 
holomorphic on $\Omega ^-(A)$.  This completes the proof. 
\vskip 4mm
\bf 6.\ Dropping the assumption on H\" older continuity\rm
\vskip 2mm
The map $\z\mapsto\z^\ast =1/\overline\z $ is the antiholomorphic 
reflection across 
$b\D $ which fixes $b\D$. Similarly, given $a\in\C$ and $\rho >0$, 
the map $\z\mapsto\z^\ast
= a+\rho^2/(\overline\z -\overline a)$ is the antiholomorphic
reflection across $b\D (a,\rho)$
which fixes $b\D (a,\rho)$ and maps a point on a ray emanating 
from $a$ at a distance $\gamma >0$ from $a$ to 
the point on the same ray at a distance $\rho^2/\gamma $ from $a$.  

If $f$ is a continuous function on $b\D (a,\rho)$ which is not
necessarily H\" older continuous
then the functions $f^+$ and $f^-$ defined by (1.5) are well
defined away from $b\D (a,\rho)$ but need not have
boundary values as we approach $b\D (a,\rho )$. The 
following still holds and is well known:
\vskip 2mm
\noindent \bf Lemma 6.1 \rm [Z, Vol.\ 1, p.\ 288]\ \it Let $f$ be a continuous 
function on $b\D (a,\rho)$. Define 
 $f^+$ and $f^-$ by (1.4). Then the function
$$
z\mapsto\ \ \left\{\eqalign{&f^+(z)+f^-(z^\ast)\ \ (z\in\D (a,\rho)) \cr
                     &f(z)\ \ \ \ \ \ \quad\quad\quad (z\in b\D (a,\rho)) 
					 \cr}\right\}
$$	
is continuous on $\overline\D (a,\rho)$. In fact, on 
$\D (a,\rho)$ it coincides with the 
Poisson integral of $f$. \rm
\vskip 2mm
\noindent We want to show that a generalization of this 
holds for continuous functions on the annulus $A$ with 
zero means on circles surrounding the origin.  

Suppose that $f$ is a continuous function on $A$ which 
satisfies (1.5) whenever $b\D (a,\rho )
\subset A$ surrounds the origin. Recall that for each 
nonnegative integer $n$ there are homogeneous 
polynomials $P_n$ and $Q_n$ of degree $n$ such that putting $z=re^{i\theta}$ into
$$
(1/\overline z)\sum_{n=0}^\infty P_n(z,1/\overline z) +
(1/ z)\sum_{n=0}^\infty Q_n(\overline z,1/z) 
$$
we get the Fourier series of the function 
$e^{i\theta}\mapsto f(re^{i\theta}),\ r_1\leq r\leq r_2$. 
We now show that one can define a holomorphic function $F^+$ on $\Omega ^+(A)$ by
$$
F^+(z,w)= (1/w)\sum_{n=0}^\infty P_n(z,1/w) \ \ \ ((z,w)\in\Omega ^+ (A))
\eqno (6.1)
$$
where the series converges uniformly on compacta in $\Omega ^+ (A)$
and a holomorphic function $F^-$ on $\Omega ^-(A)$ by
$$
F^-(z,w)= (1/z)\sum_{n=0}^\infty Q_n(w,1/z) \ \ \ ((z,w)\in\Omega ^- (A))
\eqno (6.2)
$$
where the series converges uniformly on compacta in $\Omega ^- (A)$.

Recall that for each $t,\ 0<t<1$, the series 
$(1/\overline z)\sum_{j=1}^\infty t^{j+1}P_j(z,1/\overline z)$ 
converges uniformly on $A$ which, by a reasoning similar to the 
one in Section 5 implies that the function
$$
F_t^+(z,w)=(1/w)\sum_{j=0}^\infty t^{j+1}P_j(z,1/w)
\eqno (6.3)
$$
is well defined, bounded and continuous on $\Omega^+(A)\cup b\Omega ^+(A)$
and holomorphic on $\Omega ^+(A)$ since the series (6.3) converges uniformly 
on $\Omega ^+(A)\cup b\Omega ^+ (A)$. We have $F_t^+(z,\overline z)= f_t^+(z)\ 
(z\in A)$. Similarly, the function 
$$
F^-_t(z,w)= (1/z)\sum_{j=0}^\infty t^{j+1}Q_j(w,1/z) 
$$
is well defined, bounded and continuous on $\Omega^-(A)\cup b\Omega ^-(A)$
and holomorphic on $\Omega ^-(A)$ since the series  converges uniformly 
on $\Omega ^-(A)\cup b\Omega ^- (A)$. We have $F_t^-(z,\overline z)= f_t^-(z)\ 
(z\in A)$.

Using the homogeneity of $P_j$ we rewrite (6.3) to 
$$
F_t^+(z,w)= {1\over {(w/t)}}\sum _{j=0}^\infty  P_j(tz,1/(w/t))
$$
and we see that (6.1) converges uniformly on $T_t(\Omega ^+ (A))$ where  $T_t(z,w)=(tz,w/t)$. Given 
a compact set $K\subset \Omega ^+ (A)$ there is a $t,\ 0<t<1,$ such that $K\subset 
T_t (\Omega ^+ (A))$ so (6.1) converges uniformly on compacta in $\Omega ^+ (A)$. We have
$$ 
F_t^+(z,w) = F(T_t (z,w))\ \ ((z,w)\in T_t^{-1}(\Omega ^+ (A))).
$$
Since $T_t$ converges uniformly to the identity as $t\rightarrow 1$ it follows that
 $$
F^+(z,w)= \lim_{t\rightarrow 1}F_t^+(z,w)\ \ \ \ ((z,w)\in\Omega ^+ (A)).
\eqno (6.4)
$$
where the convergence is uniform on compact sets in $\Omega ^+ (A)$. In the same way we see that 
 $$
F^-(z,w)= \lim_{t\rightarrow 1}F_t^-(z,w)\ \ \ \ ((z,w)\in\Omega ^- (A)).
\eqno (6.5)
$$
where the convergence is uniform on compact sets in $\Omega ^- (A)$.

One can verify that for each $r,\ r_1\leq r\leq r_2$, and for each $t,\ 0<t<1$, 
$$
f_t(re^{i\theta})\equiv\cP _r(te^{i\theta})\ \ (\theta\in\R)
$$
where $\cP _r$ is the Poisson integral of the function 
$e^{i\theta}\mapsto f(re^{i\theta})$.  Now,
$\cP _r(te^{i\theta })\rightarrow f(re^{i\theta })$ uniformly in $\theta $ as 
$t\nearrow 1$ 
and since $f$ is uniformly continuous on 
$A$, the standard proof of the boundary continuity of the Poisson integral shows that 
the convergence is uniform also in $r,\ r_1\leq r\leq r_2$. 
 Thus,
$$
f_t\rightarrow f \hbox{\ uniformly on \ } A \hbox{\ as \ } t\nearrow 1.
\eqno (6.6)$$ 
\vskip 4mm
\bf 7.\ Continuous functions with zero means on circles \rm
\vskip 2mm
\vskip 2mm
\noindent\bf Theorem 7.1\ \it Let $f$ be a continuous 
function on $A$ which satisfies (1.5) for each $b\D (a,\rho )\subset A$ 
which surrounds the origin. There are a holomorphic 
function $F^+$ on $\Omega ^+ (A)$ and a 
holomorphic function $F^-$ on $\Omega ^- (A)$ such that the function
$$
(z,w)\mapsto \left\{\eqalign{&F^+(z,w)+
F^-(\overline w,\overline z)\ \ ((z,w)\in\Omega^+(A))\cr
&f(z)\quad\quad\quad\quad \quad\quad\quad\quad(z,\overline z)\in\tilde A)\cr}\right\}
\eqno (7.1)
$$
has a bounded continuous extension to $\Omega ^+(A)\cup b\Omega ^+(A)$.  
\vskip 2mm\rm
\noindent Thus, for each $z\in A$ we have
$$
f(z)=\lim_{(\xi,\eta)\rightarrow (z,\overline z), (\xi,\eta)\in \Omega ^+(A)}
[F^+(\xi,\eta)+F^-(\overline\eta ,\overline\xi)] .
$$
Theorem 7.1 is the analogue of Lemma 6.1 where $b\D (a,\rho)$ 
is replaced by $\tilde A$ and $\D (a,\rho)$ is 
replaced by $\Omega ^+ (A)$. The function (7.1) is the analogue of the Poisson integral of $f$. 
It is the bounded continuous extension of $F(z,\overline z)= f(z)$ to
$\Omega ^+(A)\cup b\Omega ^+ (A)$ which is pluriharmonic on $\Omega ^+ (A)$. 
\vskip 1mm
\noindent\bf Proof of Theorem 7.1.\ \rm For each $t,\ 0<t<1$, define
$$
\Psi _t(z,w) = F_t^+(z,w)+F_t^-(\overline w,\overline z)\ \ 
\ ((z,w)\in\Omega^+(A)\cup b\Omega ^+(A)),
$$
where $F_t^+$ and $F_t^-$ are as in Section 6. 
The properties of $F_t^+$ and $F_t^-$ imply that $\Psi _t$ is bounded 
and continuous on $\Omega ^+(A)\cup b\Omega ^+(A)$, pluriharmonic on
$\Omega ^+(A)$ and satisfies $\Psi _t(z,\overline z)=f_t^+(z)+f_t^-(z) = f_t(z)\ 
(z\in A)$. It follows that for each $t, s,\ 0<t<1,\ 0<s<1$ the function
$$
(z,w)\mapsto |\Psi _t(z,w)-\Psi _s(z,w)|,
\eqno (7.2)
$$
restricted to $\Lambda_{a,\rho}\cup b\Lambda _{a,\rho}$, attains 
its maximum on $b\Lambda _{a,\rho}$ whenever $b\D (a,\rho)\subset A$ 
surrounds the origin. This is so since $\Lambda _{a,\rho}$ is 
 biholomorphically equivalent to the punctured disc and since isolated
 singularities are removable for bounded harmonic functions. Thus, the function 
 (7.2) attains its maximum on $\Omega^+(A)\cup b\Omega^+(A)$ on $\tilde A$. By (6.6) 
 the restrictions of functions $\Psi _t$ to $\tilde A$ converge uniformly as $t\nearrow 1$
 which, by the preceding discussion implies that as $t\nearrow 1$ the functions 
 $\Psi _t$ converge 
 uniformly on $\Omega^+ (A)\cup b\Omega^+(A)$ to a bounded continuous function $\Psi $ which 
 is pluriharmonic on $\Omega ^+(A)$ and which satisfies $\Psi (z,\overline z) = f(z)\ 
 (z\in A)$. Now (6.4) and (6.5) imply that $\Psi (z,w) = F^+(z,w)+F^-(z,w)\  ((z,w)\in 
 \Omega ^+(A))$ where $F^+$ and $F^-$ are given by (6.1) and (6.2). This completes the proof.
\vskip 5mm
\noindent \bf Acknowledgement\ \ \rm This work was supported in part by a grant from the
Ministry of Education, Science and Sport of the Republic of Slovenia. 
\vfill
\eject
\centerline{\bf References}
\vskip 3mm
\noindent [AG]\ M.\ Agranovsky,\  J.\ Globevnik:\  Analyticity on circles for rational
and real analytic functions of two real variables.

\noindent To appear in J.\ d'Analyse Math.
\vskip 2mm
\noindent [EK]\ C.\ L.\ Epstein,\ B.\ Kleiner:\  Spherical means in annular regions.

\noindent Comm.\ Pure Appl.\ Math.\ 46 (1993) 441-451
\vskip 2mm
\noindent [G1]\ J.\ Globevnik:\ Zero integrals on 
circles and characterizations of harmonic and analytic functions.

\noindent Trans.\ Amer.\ Math.\ Soc.\ 317 (1990) 313-330
\vskip 2mm
\noindent [G2]\ J.\ Globevnik:\ Holomorphic extensions from open families of circles.

\noindent Trans.\ Amer.\ Math.\ Soc.\ 355 (2003) 1921-1931
\vskip 2mm
\noindent [M]\ N.\ I.\ Muskhelishvili: \it Singular integral equations. \rm

\noindent Noordhoff, Groningen 1959
\vskip 2mm
\noindent [V]\ V.\ Volchkov:\ Spherical means on Euclidean spaces.

\noindent Ukrain.\ Math.\ J.\ 50 (1998) 1310-1315
\vskip 2mm
\noindent [Z]\ A.\ Zygmund:\ \it Trigonometric series.\rm 

\noindent Cambr.\ Univ.\ Press 1959
\vskip 20mm
\noindent Institute of Mathematics, Physics and Mechanics

\noindent University of Ljubljana, Ljubljana, Slovenia

\noindent josip.globevnik@fmf.uni-lj.si

\bye